\documentclass[letterpaper, 10 pt, conference]{ieeeconf}  

\IEEEoverridecommandlockouts                              

\usepackage{graphicx,array,mathrsfs,amssymb,caption,mathtools,algorithm2e} 
\usepackage[table]{xcolor}
\usepackage{amsmath} 
\newcolumntype{P}[1]{>{\centering\arraybackslash}p{#1}}
\newcolumntype{M}[1]{>{\centering\arraybackslash}m{#1}}

\title{\LARGE \bf
Clustering with Capacity and Size Constraints: A Deterministic Approach
}

\author{Mayank Baranwal$^{1}$ and Srinivasa M. Salapaka$^{2}$
\thanks{$^{1}$Mayank Baranwal is with the Department of Mechanical Engineering, University of Illinois, Urbana-Champaign, USA
        {\tt\small baranwa2@illinois.edu}}%
\thanks{$^{2}$Srinivasa M. Salapaka is with the Faculty of Mechanical Engineering, University of Illinois, Urbana-Champaign, USA
        {\tt\small salapaka@illinois.edu}}%
}

\begin{document}

\maketitle
\thispagestyle{empty}
\pagestyle{empty}

\begin{abstract}

This paper discusses a deterministic clustering approach to capacitated resource allocation problems. In particular, the Deterministic Annealing (DA) algorithm from the data-compression literature, which bears a distinct analogy to the phase transformation under annealing process in statistical physics, is adapted to address problems pertaining to clustering with several forms of size constraints. These constraints are addressed through appropriate modifications of the basic DA formulation by judiciously adjusting the free-energy function in the DA algorithm. At a given value of the annealing parameter, the iterations of the DA algorithm are of the form of a Descent Method, which motivate scaling principles for faster convergence.

\end{abstract}

\section{INTRODUCTION}\label{sec:Intro}
Clustering is an integral task to many facility location problems (FLPs) which come up in various forms in seemingly unrelated areas of coarse control quantization \cite{elia2001stabilization}, minimum distortion problem in data compression \cite{gersho2012vector}, pattern recognition \cite{therrien1989decision}, image segmentation \cite{oliver2006improving}, dynamic coverage \cite{sharma2012entropy}, neural networks \cite{haykin2004comprehensive}, graph aggregation \cite{xu2014aggregation}, and coverage control \cite{xu2010dynamic, cortes2002coverage}. These problems each with different and unrelated goals, in fact have some fundamental common attributes - (1) obtaining an optimal partition of the underlying domain, and (2) optimally assigning values from a finite (or countable) set to each cell of the partition. The differences in these problems come from having different conditions of optimality and constraints. For example, the image segmentation problem consists of optimally partitioning an image into connected components. Similarly, the coarse control quantization problem requires partitioning the state-space and allocating a control value to each partition such that the stability of the underlying dynamical system is guaranteed.

These optimization problems are largely non convex, computationally complex and suffer from multiple poor local minima that riddle the cost surface \cite{gray1982multiple}. Many heuristics have been proposed to address these difficulties. These include repeated optimization with random initialization (such as Lloyd's algorithm \cite{lloyd1982least}) and efficient rules for cluster splits and merges. In this context, {\em simulated annealing} (SA) algorithm \cite{kirkpatrick1983optimization}, which is motivated by an analogy to the statistical mechanics of annealing (and phase transformations) in solids, was shown to be an effective iterative metaheuristic probabilistic technique for approximating the {\em global} optimum of an optimization problem, however with an annealing schedule so slow that the algorithm loses practicality in many applications.

In this paper, we discuss the Deterministic Annealing (DA) algorithm developed in the data-compression literature \cite{rose1990deterministic, rose1998deterministic}. The DA algorithm enjoys the best of both the worlds. On one hand it is {\em deterministic}, i.e., it does not wander randomly on the energy surface. On the other hand, it is still an annealing method designed to aim at the global optimum, instead of going directly to a near local minima. The DA algorithm, thus, has an ability to avoid poor local optima while still maintaining a relatively faster annealing schedule. It formulates an effective {\em free-energy} function parameterized by an annealing parameter and this function is deterministically optimized at successively increased values of the annealing parameter.

While the original DA algorithm was developed in the context of unconstrained clustering (i.e., clustering without any capacity constraint), most real-world problems demand limited capacities on the cluster (resource) sizes, such as limited link capacity in communication channels, and bounded capacities of vehicles in pick-up and delivery problems (PDPs). Earlier approaches using DA have attempted modifications to accommodate capacity constraints heuristically \cite{salapaka2003constraints, sharma2008entropy}, which typically do not quantitatively satisfy the capacity  constraints - especially the cluster-size constraints. In this paper, we describe the modifications that we have made to adapt the DA algorithm to address several heterogeneous capacity constraints. These modifications are designed such that solutions do satisfy these constraints; we substantiate their effectiveness through some practical scenarios. We also derive some spatial scaling laws to obtain faster convergence rate for the DA algorithm, thereby making the algorithm tractable for handling complex combinatorial problems.

The rest of the paper is structured as follows. Section \ref{sec:problem_form} introduces the clustering problem and its variants involving constraints on cluster sizes. This is then followed by an introduction the Deterministic Annealing framework in Section \ref{sec:DA_intro}. We then discuss several modifications of the DA algorithm to address capacity constraints of multiple types, which is followed by application of the capacity constrained DA on some practical scenarios. Section \ref{sec:scaling} discusses scaling laws and convergence rate associated with the DA algorithm, followed by the section on conclusions and future work.

{\it Remark}: For any two natural numbers $a, b\in\mathbb{N}$, the index set $\{a,a+1,\hdots,b\}$ is denoted by $\mathbb{N}_{[a,b]}$. The set of non-negative real numbers is denoted by $\mathbb{R}_+$. The distance between two vectors $x,y\in\mathbb{R}^n$ is considered to be squared-Euclidean unless stated otherwise.], i.e., $d(x,y) = \|x-y\|_2^2$.

\section{PROBLEM FORMULATION}\label{sec:problem_form}
The problems discussed in this paper have the objective of obtaining resource (cluster) locations $\mathcal{Y}=\{y_j:j\in\mathbb{N}_{[1,K]}\}$ for a given set of demand points located at $\mathcal{X}=\{x_i:i\in\mathbb{N}_{[1,N]}\}$ such that the total sum of squared distances from each demand point to its nearest resource is minimum. The demand point follow the probability distribution $p(x_i), x_i\in\mathcal{X}$. This results in the following classes of optimization problems:

{\bf P1. No constraints:}
In this setting, the objective is to find $K$ resource locations such that the cumulative distance from each demand point to its nearest resource is minimized; i.e., 
\begin{equation}\label{eq:P1}
	\min\limits_{y_j,j\in\mathbb{N}_{[1,K]}}\sum\limits_{i\in\mathbb{N}_{[1,N]}}p(x_i)\left\{\min\limits_{j\in\mathbb{N}_{[1,K]}}\|x_i-y_j\|_2^2\right\}
\end{equation}
The classical DA algorithm by Rose \cite{rose1998deterministic} directly solves this class of problems. This formulation does not take into account any constraints on the {\em size} and {\em type} of resources. Equivalently the problem is viewed as obtaining the partition of the domain $\mathcal{X}$ of demand points into $K$ clusters $C_j\subset\mathcal{X}, \forall j\in\mathbb{N}_{[1,K]}$ such that $\cup_jC_j = \mathcal{X}; C_{j_1}\cap C_{j_2} = \phi, \forall j_1\neq j_2$, and allocate a representative resource to each cluster such that the expected distance of the demand point $x_i\in C_j$ from the corresponding resource location $y_j$ is minimized, i.e.
\begin{equation}\label{eq:alternate}
	\min\limits_{\{y_j\}, \{C_j\}}\sum\limits_{j\in\mathbb{N}_{[1,K]}}\sum\limits_{x_i\in C_j}p(x_i)\|x_i-y_j\|_2^2.
\end{equation}
Note that these optimization problems are non convex and computationally complex. For instance, the optimal allocation of $20$ resources in a domain of demand $30$ points would require search over $30$ million partitions, which renders most deterministic optimization approaches practically infeasible. The constrained problems result in following formulations:\\

{\bf P2. Heterogeneous capacity constraints:}\\
In this scenario, the objective is to obtain a partition of the underlying domain of demand points, such that the {\em relative size} of the $j^{th}$ cluster $C_j$ is equal to a pre-specified value $\lambda_j$. Moreover, we assume $\sum\limits_{j\in\mathbb{N}_{[1,K]}}\lambda_j = 1$.
\begin{eqnarray}\label{eq:P2}
	&&\min\limits_{y_j,j\in\mathbb{N}_{[1,K]}}\sum\limits_{i\in\mathbb{N}_{[1,N]}}p(x_i)\left\{\min\limits_{j\in\mathbb{N}_{[1,K]}}\|x_i-y_j\|_2^2\right\} \nonumber \\
	&&\text{s.t.} \qquad \quad \sum\limits_{i\in\mathbb{N}_{[1,N]}}p(x_i)v_{ij} = \lambda_j, \quad j\in\mathbb{N}_{[1,K]} \\
	&&\text{where,} \qquad \quad v_{ij} = \left\{
		\begin{array}{cc}
			1, & \text{if} \quad x_i\in C_j\\
			0, & \text{else}
		\end{array}
	\right. \nonumber
\end{eqnarray}
Note that the constraints here seem to be independent of the cost function. This is not so since the cluster size $\lambda_j$ depends on the partition $\{C_j\}$, which is a decision variable for the optimization problem (\ref{eq:P2}). This parameter is more naturally incorporated in the cost function in the modified algorithm that is presented in Section \ref{sec:DA_cap}.\\

{\bf P3. Capacity constraints with multiple types of demand points:}\\
These constraints relate to the scenario where the demand points are heterogeneous and have multiple types $k\in\mathbb{N}_{[1,p]}$ and the resources have capacity constraints given by $\lambda_{jk}$. Here $\lambda_{jk}$ denotes the capacity of resource $j$ for demand point of type $k$. We have the following optimization problem to be solved:
\begin{equation}\label{eq:P3}
	\begin{aligned}
		&\min\limits_{y_j,j\in\mathbb{N}_{[1,K]}}\sum\limits_{i\in\mathbb{N}_{[1,N]}}p(x_i)\left\{\min\limits_{j\in\mathbb{N}_{[1,K]}}\|x_i-y_j\|_2^2\right\} \\
		&\text{s.t.} \quad \sum\limits_{i\in\mathbb{N}_{[1,N]}}p(x_i)v_{ijk} = \lambda_{jk}, \quad j\in\mathbb{N}_{[1,K]}, k\in\mathbb{N}_{[1,p]}\\
		&\text{where,} \qquad v_{ijk} = \left\{
		\begin{array}{cc}
			1, & \text{if} \quad \text{$x_i\in C_j$ and $x_i$ is of type-$k$} \\
			0, & \text{else}
		\end{array}
	\right.
	\end{aligned}
\end{equation}
Note that the capacitated problems inherit the complexity of P1. In this paper, we describe an iterative algorithm (developed in \cite{rose1990deterministic}) and its modifications to address various capacity constraints, while still able to obtain {\em good} solutions. 

\section{DETERMINISTIC ANNEALING ALGORITHM: A MAXIMUM-ENTROPY PRINCIPLE APPROACH FOR CLUSTERING}\label{sec:DA_intro}
At its core, the Deterministic Annealing (DA) algorithm solves a facility location problem (FLP), i.e., given a set of demand point locations $\mathcal{X}=\{x_i,i\in\mathbb{N}_{[1,N]}\}$, find $K\in\mathbb{N}$ facility locations $\mathcal{Y}=\{y_j,j\in\mathbb{N}_{[1,K]}\}$ such that the {\em total weighted sum of the distance of each demand point from its nearest facility location is minimized}. The FLP is mathematically described in (\ref{eq:P1}). Borrowing from the data compression literature \cite{cover2012elements}, we define {\em distortion} as a measure of the average distance of a demand point to its nearest facility, given by $D(\mathcal{X},\mathcal{Y}) = \sum\limits_{i\in\mathbb{N}_{[1,N]}}p(x_i)\min\limits_{j\in\mathbb{N}_{[1,K]}}d(x_i,y_j)$. The solution to an FLP satisfies the following two necessary (but not necessarily sufficient) properties:
\begin{itemize}
	\item Voronoi partitions: The partition of the domain is such that each demand point in the domain is associated only to its nearest resource (cluster) location.
	\item Centroid condition: The resource location $y_j$ is at the centroid of the $j^{th}$ cluster $C_j$.
\end{itemize}

Most algorithms for FLP (such as Lloyd's \cite{lloyd1982least}) are overly sensitive to the initial resource locations. This is primarily due to the distributed aspect of the FLPs, where any change in the location of the $i^{th}$ demand point affects $d(x_i,y_j)$ only with respect to the {\em nearest} facility $j$. The DA algorithm suggested by Rose \cite{rose1998deterministic}, overcomes this sensitivity by allowing {\em fuzzy} association of every demand point to each facility through an association probability, $p(y_j|x_i)$:  
\begin{equation}\label{eq:mod_D}
	\bar{D}(\mathcal{X},\mathcal{Y}) = \sum\limits_{i\in\mathbb{N}_{[1,N]}}p(x_i)\sum\limits_{j\in\mathbb{N}_{[1,K]}}p(y_j|x_i)d(x_i,y_j).
\end{equation}
Thus the notion of average distance of a demand point from its nearest facility is replaced by the weighted average distance of demand points to {\em all} the facilities. The probability distribution $\{p(y_j|x_i)\}$ determines the trade-off between decreasing the {\em local} influence and the deviation of the modified distortion $\bar{D}$ from the original distortion measure $D$. The uncertainties in facility locations $\{y_j\}$ with respect to the demand point locations $\{x_i\}$ is captured by Shannon {\em entropy} $H(\mathcal{Y}|\mathcal{X}) = -\sum\limits_{i\in\mathbb{N}_{[1,N]}}p(x_i)\sum\limits_{j\in\mathbb{N}_{[1,K]}}p(y_j|x_i)\log(p(y_j|x_i))$, widely used in data compression literature \cite{cover2012elements}. Therefore, maximizing the entropy is commensurate with decreasing the {\em local} influence.

This trade-off between maximizing the entropy and minimizing the distortion in Eq. (\ref{eq:mod_D}) is addressed by seeking the probability distribution $\{p(y_j|x_i)\}$ that minimize the {\em free-energy}, or the Lagrangian, given by $F\coloneqq\bar{D}(\mathcal{X},\mathcal{Y})-\dfrac{1}{\beta}H(\mathcal{Y}|\mathcal{X})$, where $\beta$ is the Lagrange multiplier and bears a direct analogy to the inverse of the {\em temperature} variable in an annealing process. The association weights $\{p(y_j|x_i)\}$ that minimize the free-energy function are given by the {\em Gibbs} distribution
\begin{equation}\label{eq:Gibbs}
	p(y_j|x_i) = \frac{e^{-\beta d(x_i,y_j)}}{\sum\limits_{j\in\mathbb{N}_{[1,K]}}e^{-\beta d(x_i,y_j)}}.
\end{equation}
By substituting the Gibbs distribution (\ref{eq:Gibbs}), the corresponding {\em free-energy} function is obtained as
\begin{equation}\label{eq:Free_energy}
	F(\mathcal{Y}) = -\frac{1}{\beta}\sum\limits_{i\in\mathbb{N}_{[1,N]}}p(x_i)\log\bigg(\sum\limits_{j\in\mathbb{N}_{[1,K]}}e^{-\beta d(x_i,y_j)}\bigg).
\end{equation}
In the DA algorithm, the free-energy function is {\em deterministically} optimized at successively increased values of the annealing parameter $\beta$.

The readers are encouraged to refer to \cite{parekh2015deterministic} for detailed analysis on the complexity of the DA algorithm. For implementation on very large datasets, a scalable modification of the DA is proposed in \cite{sharma2006scalable}.

\section{CAPACITY CONSTRAINED DA}\label{sec:DA_cap}
In this section, we develop an iterative algorithm to address capacity constraints discussed in P2-P3 in Section \ref{sec:problem_form}. Adaptation of the DA algorithm for handling capacity constraints in the context of locational optimization problems was earlier discussed in \cite{salapaka2003constraints}. However, the approach adopted in \cite{salapaka2003constraints} satisfies the constraints only under the assumption of the uniform distribution of the demand points, and thus renders the algorithm impotent to majority of the real-world problems with non-uniform distribution of demand points. This approach is discussed further in Section \ref{sec:sim_results}, where the approach adopted in \cite{salapaka2003constraints} results in clusters with highly inconsistent constraint-matching. We now describe  the methodology adopted in this paper for addressing these constraints.

{\bf P2. Heterogeneous capacity constraints:}\\
In the setting, we consider resources $j\in\mathbb{N}_{[1,K]}$ with relative heterogeneous capacities $\lambda_j\in[0,1]$ and the objective is to allocate resources $y_j$ to the demand points $x_i\in\mathcal{X}$. We address this in the DA framework by requiring
\begin{equation}\label{eq:sol1_1}
	p(y_1):p(y_2):\hdots:p(y_K) = \lambda_1:\lambda_2:\hdots:\lambda_K
\end{equation}
where $p(y_j) = \sum\limits_{i\in\mathbb{N}_{[1,N]}}p(x_i)p(y_j|x_i)$ is the {\em mass} associated with cluster $C_j$. The capacity constraints are incorporated into the DA framework through a modified Gibbs distribution given by
\begin{equation}\label{eq:sol1_2}
	p(y_j|x_i) = \frac{\eta_je^{-\beta d(x_i,y_j)}}{\sum\limits_{j\in\mathbb{N}_{[1,K]}}\eta_je^{-\beta d(x_i,y_j)}}
\end{equation}
where $\eta_j\in[0,1]\forall j\in\mathbb{N}_{[1,K]}$ specifies the relative weight of the $j^{th}$ resource $y_j$ and can be interpreted as the number of copies of $y_j$ in the cluster $C_j$. During {\em fuzzy} initialization when each demand point is uniformly associated with every resource (i.e., $\beta\approx 0$), $\eta_j$ is initialized to $\lambda_j$ and therefore $p(y_j) = \lambda_j$ in the beginning of the annealing process.

The {\em free-energy} function is modified as
\begin{equation}\label{eq:sol1_3}
	F(\mathcal{Y},\eta) = -\frac{1}{\beta}\sum\limits_{i\in\mathbb{N}_{[1,N]}}p(x_i)\log\bigg(\sum\limits_{j\in\mathbb{N}_{[1,K]}}\eta_je^{-\beta d(x_i,y_j)}\bigg).
\end{equation}
The update equation for the resource location $y_j, j\in\mathbb{N}_{[1,K]}$ is obtained by setting the derivative of the modified free-energy function w.r.t. $y_j$ to zero, which results in
\begin{equation}\label{eq:sol1_4}
	y_j = \frac{\sum\limits_{i\in\mathbb{N}_{[1,N]}}p(x_i)p(y_j|x_i)x_i}{p(y_j)}.
\end{equation}
Note that the update equations for the resource locations have implicit dependence on the weight parameters $\eta_j$, which in turn are again coupled with resource locations $y_j$ and annealing parameter $\beta$ through cluster probabilities $p(y_j)$. The cluster probability (weight) $p(y_j)$ is given by
\begin{eqnarray}\label{eq:sol1_5}
	p(y_j) &=& \sum\limits_{i\in\mathbb{N}_{[1,N]}}p(x_i)p(y_j|x_i) \nonumber\\
	&=& \sum\limits_{i\in\mathbb{N}_{[1,N]}}p(x_i)\frac{\eta_je^{-\beta d(x_i,y_j)}}{\sum\limits_{j\in\mathbb{N}_{[1,K]}}\eta_je^{-\beta d(x_i,y_j)}}\nonumber\\
	\Rightarrow \eta_j &=& \frac{p(y_j)}{\sum\limits_{i\in\mathbb{N}_{[1,N]}}p(x_i)\frac{e^{-\beta d(x_i,y_j)}}{\sum\limits_{j\in\mathbb{N}_{[1,K]}}\eta_je^{-\beta d(x_i,y_j)}}}
\end{eqnarray}
Noting that the desired {\em mass} associated with the resource $y_j$ is $\lambda_j$, i.e., $p(y_j)=\lambda_j, \forall j\in\mathbb{N}_{[1,K]}$, the update equation for $\eta_j$ is obtained as
\begin{equation}\label{eq:sol1_6}
	\eta_j = \frac{\lambda_j}{\sum\limits_{i\in\mathbb{N}_{[1,N]}}p(x_i)\frac{e^{-\beta d(x_i,y_j)}}{\sum\limits_{j\in\mathbb{N}_{[1,K]}}\eta_je^{-\beta d(x_i,y_j)}}}.
\end{equation}
In the capacitated-DA algorithm, we deterministically optimize the free-energy function in (\ref{eq:sol1_3}) at successive $\beta$ values by alternating between the Eq. (\ref{eq:sol1_4}) and Eq. (\ref{eq:sol1_6}) until convergence.\\

{\bf P3. Capacity constraints with multiple types of demand points:}\\
In this setting, there are $p$-types of demand points and there are capacity-constraints on size of each cluster  for type of demand point. We use $\lambda_{jk},\forall j\in\mathbb{N}_{[1,K]}, k\in\mathbb{N}_{[1,p]}$ to denote the capacity requirements on resource $j$ for demand point of type $k$. Similar to the previous scenario, the modified Gibbs distribution and the update equations are given by
\begin{eqnarray}\label{eq:sol2_1}
	p(y_j|x_i) = \frac{\sum\limits_{k\in\mathbb{N}_{[1,p]}}\eta_{jk}e^{-\beta d(x_i,y_j)}}{\sum\limits_{j\in\mathbb{N}_{[1,K]}}\sum\limits_{k\in\mathbb{N}_{[1,p]}}\eta_{jk}e^{-\beta d(x_i,y_j)}} \nonumber\\
	y_j = \frac{\sum\limits_{i\in\mathbb{N}_{[1,N]}}p(x_i)p(y_j|x_i)x_i}{p(y_j)} \nonumber \\
	\eta_{jk} = \frac{\lambda_{jk}}{\sum\limits_{i\in\mathbb{N}_{[1,N]}}p(x_i)\frac{e^{-\beta d(x_i,y_j)}}{\sum\limits_{j\in\mathbb{N}_{[1,K]}}\sum\limits_{k\in\mathbb{N}_{[1,p]}}\eta_{jk}e^{-\beta d(x_i,y_j)}}}.
\end{eqnarray}

\section{SIMULATION RESULTS}\label{sec:sim_results}
In this section, we consider the application of the DA algorithm and its proposed adaptations for the problem types P1-P3 through some real-world instances. In particular, we consider the following problems:

\subsection{P1. Image clustering and segmentation}\label{subsec:im_cluster}
In this scenario, we are given an image containing $N$ pixels (demand points). The {\em location} of the $i^{th}$ pixel, $x_i$, corresponds to the RGB value of the pixel, i.e., $x_i\in\mathbb{R}^3$. In image segmentation, the objective is to partition a digital image into $K$ segments (also known as superpixels) in order to simplify and change the representation of the image into something meaningful and easier to analyze. Furthermore, if each of the color components requires one byte of storage, then the size of the original image is equal to $3N$ bytes, whereas the segmented image can be represented using $3K$ bytes, thereby resulting in image compression without {\em significant} reduction in the {\em quality} of the image.

\begin{figure}
	\centering
	\includegraphics[scale=0.2]{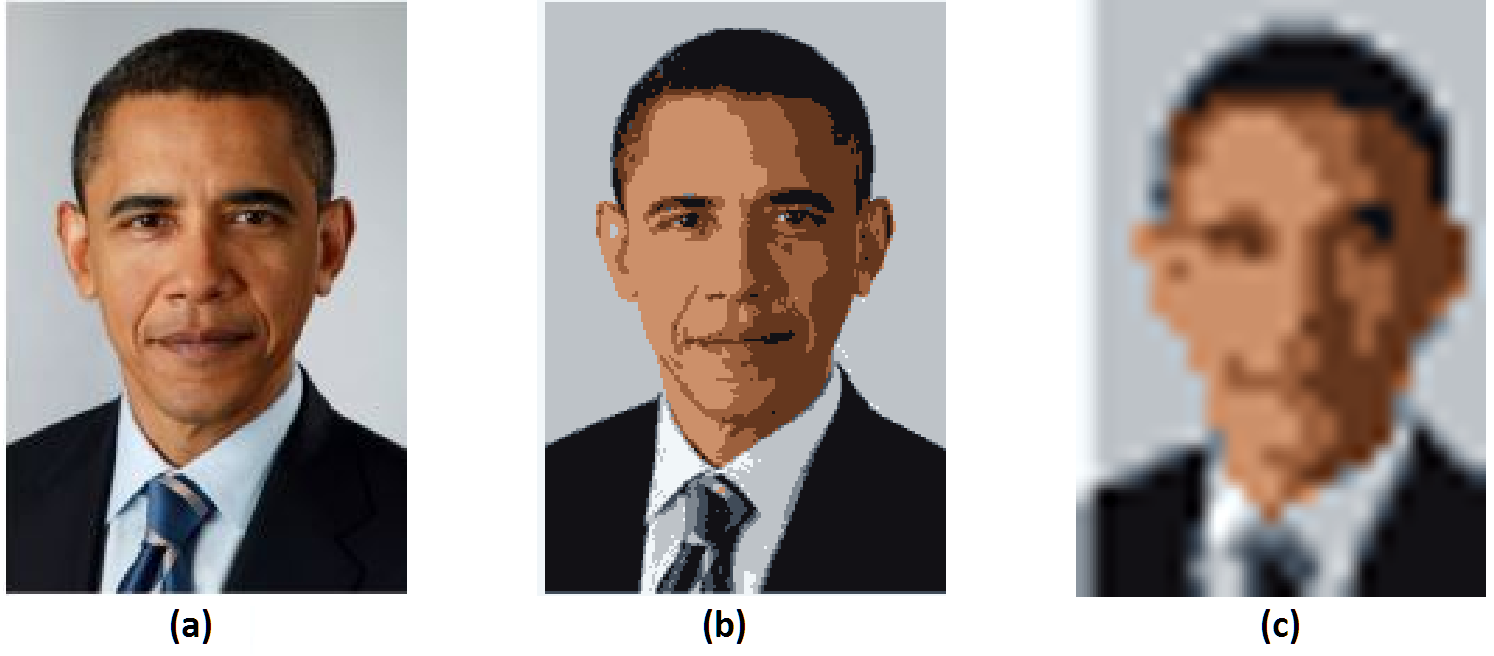}
	\caption{\small{Image segmentation using Deterministic Annealing based clustering. (a) Original image, (b) 8-bit image obtained by clustering the original image into 8 clusters, (c) Pixelated image using nearest-neighbor representation.}}
	\label{fig:obama_cluster}
	\vspace{-2em}
\end{figure}

For the simulation, we consider an image (see Fig. \ref{fig:obama_cluster}a) of dimensions $213\times 146$, i.e., $N = 31098$. We aim to represent the original image by just $K = 8$ superpixels, which results a significant reduction in the size of the image. In this context, the resource locations $\{y_j\}$ correspond to the RGB values of the superpixels. The DA formulation described in Section \ref{sec:DA_intro} directly solves the segmentation problem. The segmented image with $8$ superpixels is shown in Fig. \ref{fig:obama_cluster}b. As is seen in the segmented image, the segments (clusters) retain the {\em important} features of the original image, while reducing the size of the original image by a factor of $\sim 3887$.

Furthermore, we use the superpixels obtained using the DA approach to produce very low-resolution output image of dimensions $30\times 20$ (see Fig. \ref{fig:obama_cluster}c). Such pixel arts are often utilized by major commercial companies to convey information on compact screens and making avatars for social networks. This is achieved by averaging the RGB values over a blob of pixels in the original image and replacing the entire blob with a superpixel with the nearest RGB value. As is seen in the pixelated image in Fig. \ref{fig:obama_cluster}c, the facial features and the sharp edges are retained.

\subsection{P2. Allocating vehicles with heterogeneous capacities to service shipments at demand points}\label{subsec:cap}
In this scenario, we are given a set of $N\in\mathbb{N}$ demand points located at $\mathcal{X}=\{x_i,i\in\mathbb{N}_{[1,N]}\}$ and the objective is to cluster them into $K\in\mathbb{N}$ such that the {\em sizes} of the clusters are in the ratio $\lambda_1:\lambda_2:\hdots:\lambda_K$. This setting reflects the scheme where larger vehicles are required to serve relatively larger number of demand points.  

\begin{figure*}[!t]
	\begin{center}
	\begin{tabular}{cc}
	\includegraphics[width=0.95\columnwidth]{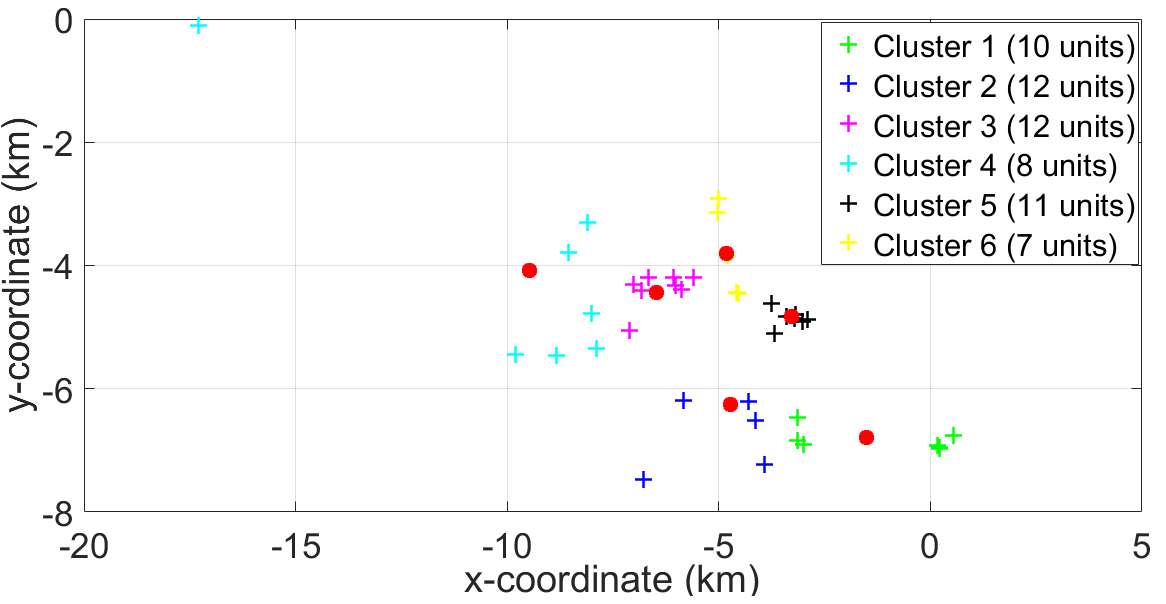}&\includegraphics[width=0.96\columnwidth]{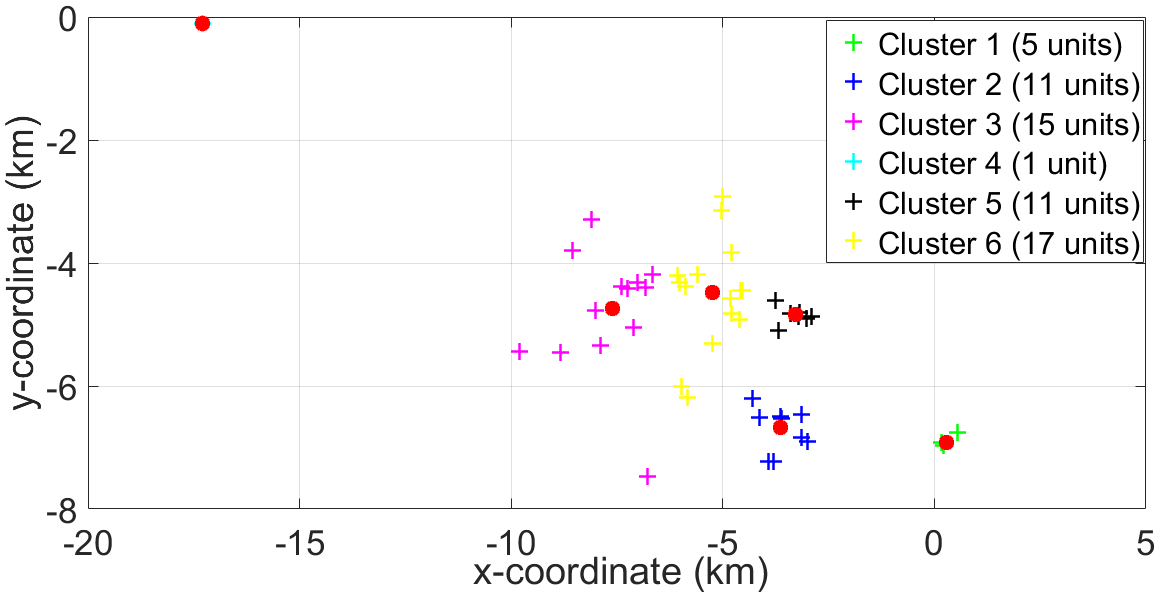}\cr
	(a)&(b)\end{tabular}
	\vspace{-0.5em}
	\caption{{\small Allocation of vehicles to customers according to pre-specified capacities. Note that some of the demand points have multiplicity more than one. (a) The capacitated-DA approach clusters the demand points in the pre-specified ratio $10:12:12:8:11:7$. (b) Results from the DA modification proposed in \cite{salapaka2003constraints}. Clearly the algorithm fails to find the clusters with desired sizes due to non-uniform distribution of demand points. In fact, the demand point in the top-left becomes a cluster in itself.}}
	\label{fig:P2}
	\vspace{-2em}
	\end{center}
\end{figure*}

For the simulation, we consider a 60 customers (demand points) data from a logistics company (name withheld). We are given a fleet of 6 vehicles with capacities of 10, 12, 12, 8, 11 and 7 units respectively. The proposed modification of the DA algorithm results in automatic allocation of vehicles to the demand points in the desired ratio. The algorithm is also compared against the previous capacitated-DA proposed in \cite{salapaka2003constraints} where $\eta_j\equiv\lambda_j,\forall j\in\mathbb{N}_{[1,K]}$ at all values of the annealing parameter, i.e., $\eta_j$ are not updated within inner-iterations of the DA algorithm. The resulting solution from the previous approach does not satisfy the capacity requirements (in fact, results in very poor quality solution due to non-uniform distribution of demand points). These results are summarized in Fig. \ref{fig:P2}. The proposed algorithm also outperforms the constrained k-means algorithm described in \cite{geetha2009improved} with smaller value of the cost function.

\subsection{P3. Pick-up problem with multi-type capacity constraints}\label{subsec:multi_cap}
In this setting, we consider a pick-up problem from a single depot and the objective is to find the arrival times of the $K\in\mathbb{N}$ vehicles (resources) to pick up $N\in\mathbb{N}$ shipments (demand points). Each shipment is equipped with a time-window $[t_{i,start},t_{i,end}]$ and belongs to one of the $p\in\mathbb{N}$ types. Additionally, the total amount of resource (vehicle) $j$ needing to be allocated to all shipments of $k^{th}$-type are $\lambda_{jk}, j\in\mathbb{N}_{[1,K]},k\in\mathbb{N}_{[1,p]}$. In this context, the {\em location} of the shipment (demand point) $x_i$ is chosen as the mid-point of the time-window, i.e., $x_i = 0.5(t_{i,start}+t_{i,end})$. $p(y_j|k)$ denotes the allocation of resource $j$ to shipment of type-$k$ and is given by
\begin{equation*}
	p(y_j|k) = \sum\limits_{i\in\mathbb{N}_{[1,N]}}p(x_i)\frac{\eta_{jk}e^{-\beta d(x_i,y_j)}}{\sum\limits_{j\in\mathbb{N}_{[1,K]}}\eta_{jk}e^{-\beta d(x_i,y_j)}}.
\end{equation*}

\begin{figure*}[!t]
	\begin{center}
	\begin{tabular}{cc}
	\includegraphics[width=0.92\columnwidth]{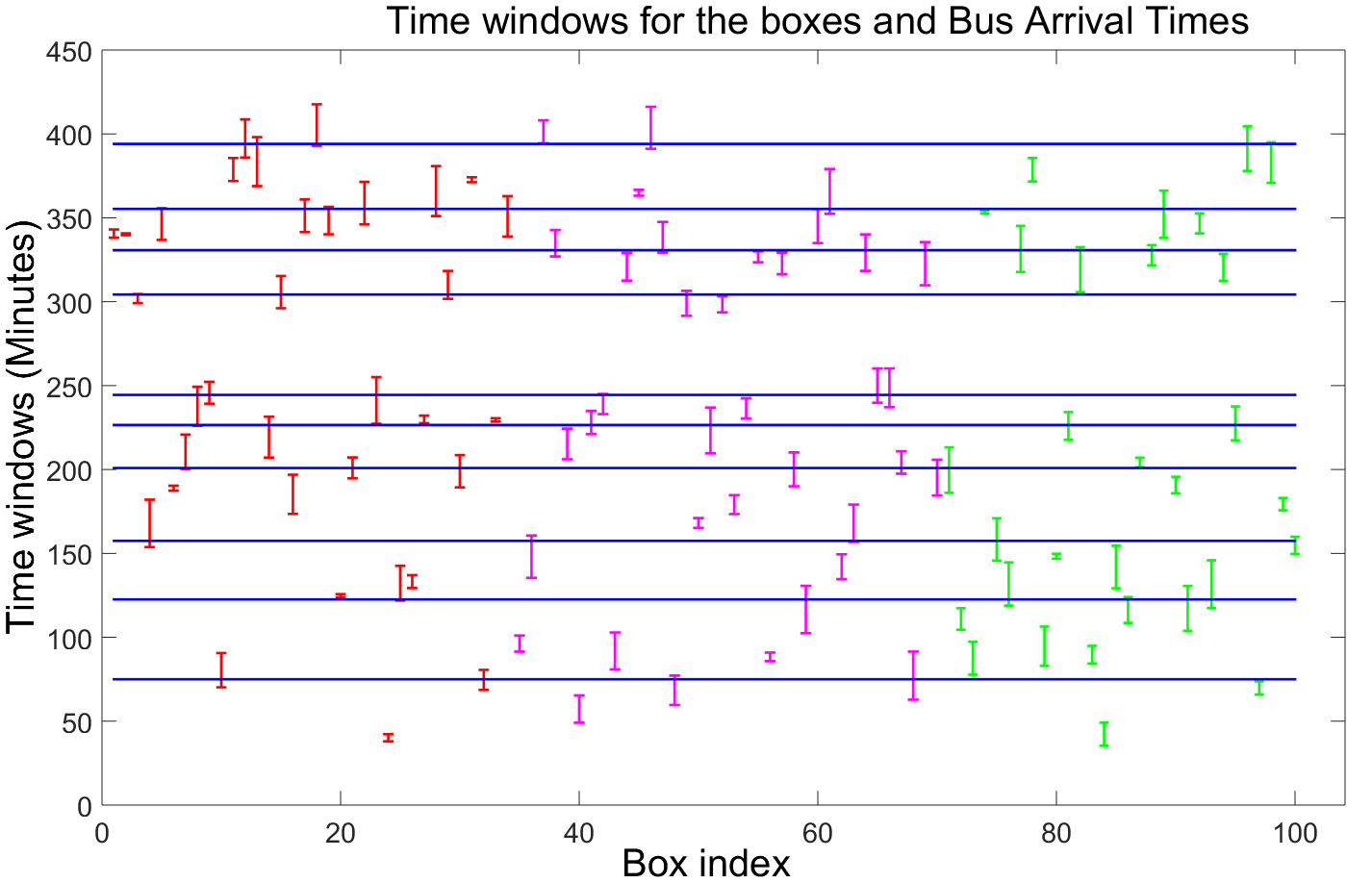}&\includegraphics[width=0.98\columnwidth]{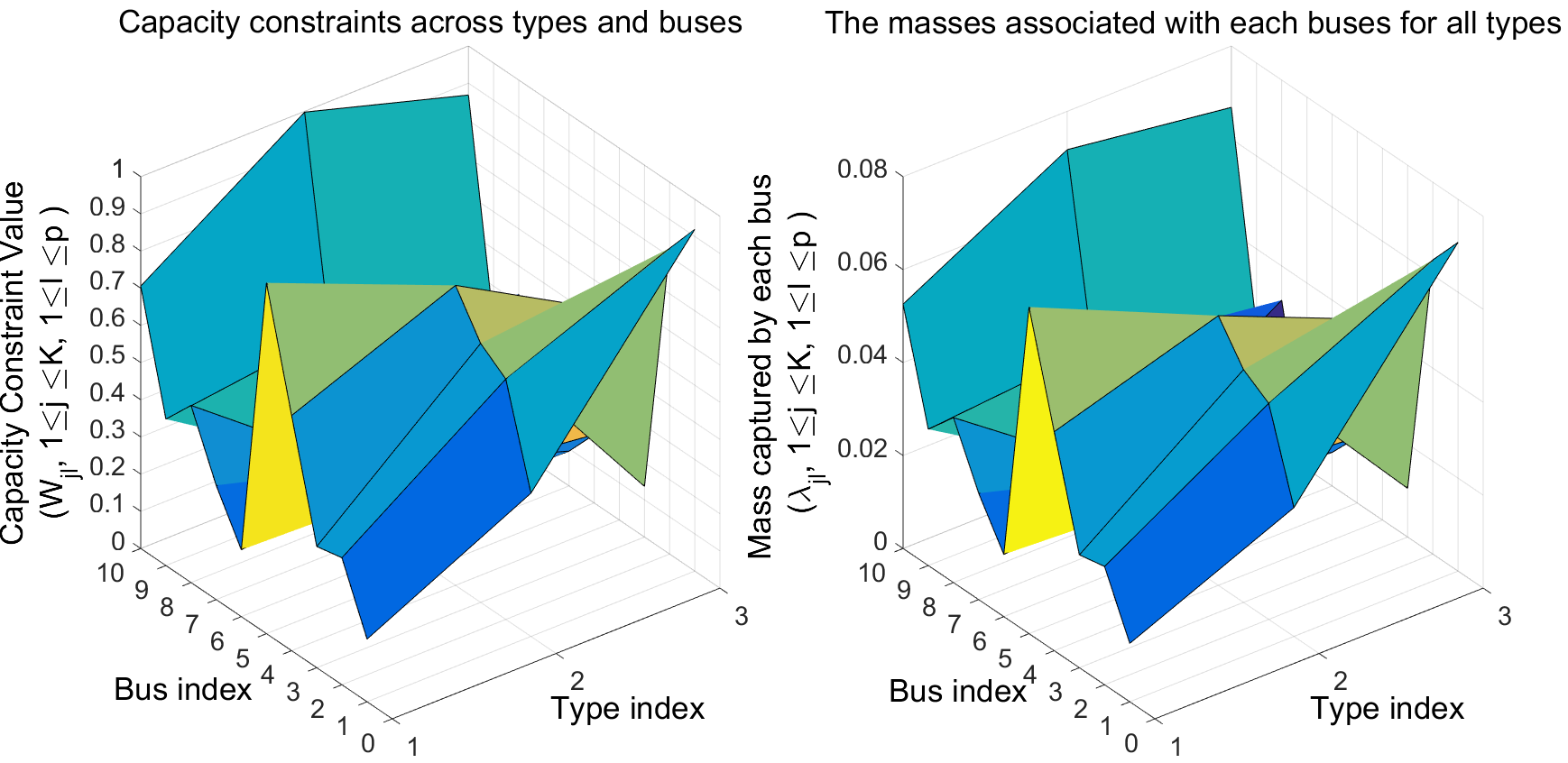}\cr
	(a)&(b)\end{tabular}
	\vspace{-0.5em}
	\caption{{\small (a) Time-Windows for 100 shipments (3 types) and the Arrival Times (Cluster locations (blue lines)) for 10 vehicles. The three colors (red, magenta and green) are used to indicate the type of shipments. (b) [Right]: Capacity constraint values for each type of shipment and vehicle. [Left]: Capacity associated with
each vehicle and each type of shipment. Clearly, the shapes indicate that the capacities obtained by the constrained DA algorithm are in proportion with the capacities provided as constraints.}}
	\label{fig:P3}
	\vspace{-2em}
	\end{center}
\end{figure*}

For the simulation, we choose an instance with 100 shipments in total, 3 types of shipments and only 10 vehicles for pickup. There are 34 shipments of type-1 (red), 36 shipments of type-2 (magenta), and 30 shipments of type-3 (green). The capacities are randomly chosen and subjected to these capacities, the algorithm works well. In Fig. \ref{fig:P3}a, the blue lines depict the vehicle arrival times, and the red, magenta and green bars indicate the timewindows for the three types of shipments respectively. Fig. \ref{fig:P3}b shows the effectiveness of the DA algorithm in respecting the capacity constraints. The plot on the left shows the capacity constraints ($\lambda_{jk}$) for each type of shipment and vehicle. The plot on the right shows the clustered mass information $p(y_j|k)$. Clearly, the shapes of the two plots match each other. Also, the numerical results find that the ratios of the capacity and clustered mass are equal for all the vehicles, i.e., for all types of shipments $k\in\mathbb{N}_{[1,p]}$, we find $p(y_j|k) = \lambda_{jk}, \forall j\in\mathbb{N}_{[1,K]}$.

\section{SCALING LAWS AND CONVERGENCE RATES}\label{sec:scaling}
We now describe the {\em spatial} scaling laws for faster convergence of the DA algorithm. In particular, we first demonstrate that the inner-iterations on the resource location $y_j$ correspond to a Descent Method and the rate of descent depends on a scaling factor $\sigma$. Utilizing the fact that the optimization problem is scale independent, we describe how by appropriately scaling the demand point locations $\{x_i\}$, the DA implementation can be expedited.

In the classical DA algorithm, at each iteration or equivalently for each value of the annealing parameter $\beta_k$, where $k$ signifies the $k^{th}$ iteration, it is required to solve the following set of implicit equations
\begin{small}
\begin{equation}\label{eq:scale_1}
	y_j = \frac{\sum\limits_{i\in\mathbb{N}_{[1,N]}}p(x_i)\frac{e^{-\beta_k\|x_i-y_j\|_2^2}}{\sum\limits_{l\in\mathbb{N}_{[1,K]}}e^{-\beta_k\|x_i-y_l\|_2^2}}x_i}{\sum\limits_{i\in\mathbb{N}_{[1,N]}}p(x_i)\frac{e^{-\beta_k\|x_i-y_j\|_2^2}}{\sum\limits_{l\in\mathbb{N}_{[1,K]}}e^{-\beta_k\|x_i-y_l\|_2^2}}}, \quad \forall j\in\mathbb{N}_{[1,K]}.	
\end{equation}
\end{small}
Eq. (\ref{eq:scale_1}) is a consequence of the {\em centroid condition}, i.e., the resource $j$ is located at the centroid of the $j^{th}$ cell. This corresponds to the following iteration scheme
\begin{small}
\begin{equation}\label{eq:scale_2}
	y_j(n+1) = \frac{\sum\limits_{i}p(x_i)p_k(y_j(n)|x_i)}{p_k(y_j(n))}\eqqcolon g_j^k(\mathcal{Y}(n)),\quad \forall j\in\mathbb{N}_{[1,K]},
\end{equation}
\end{small}
where $p_k(y_j(n)|x_i)$ is the Gibbs distribution given by $\dfrac{e^{-\beta_k\|x_i-y_j(n)\|_2^2}}{\sum\limits_{l\in\mathbb{N}_{[1,K]}}e^{-\beta_k\|x_i-y_l(n)\|_2^2}}$, $n=0,1,2,\hdots$, and $y_j(0)$ is assigned the solution of Eq. (\ref{eq:scale_1}) at the previous value of the annealing parameter $\beta_{k-1}$. Thus the above iteration scheme (\ref{eq:scale_2}) can be written as $\mathcal{Y}(n+1) = g_j^k(\mathcal{Y}(n))$. Note that the free-energy at $n^{th}$ iteration is given by $F_k(n) = -\dfrac{1}{\beta_k}\sum\limits_{i\in\mathbb{N}_{[1,N]}}p(x_i)\log\left(\sum\limits_{j\in\mathbb{N}_{[1,K]}}e^{-\beta_k\|x_i-y_j(n)\|_2^2}\right)$ and hence
\begin{eqnarray}
	\frac{1}{2}\frac{\partial F_k(n)}{\partial y_j(n)} &=& \sum\limits_{i\in\mathbb{N}_{[1,N]}}p(x_i)p_k(y_j(n)|x_i)(y_j(n)-x_i)\nonumber\\ 
	&=& p_k(y_j(n))(y_j(n) - g_j^k(\mathcal{Y}(n)))\nonumber\\
	\Rightarrow \frac{1}{2}\nabla F_k(n) &=& P_k(n)(\mathcal{Y}(n) - \mathcal{Y}(n+1))\nonumber\\
	\Rightarrow \mathcal{Y}(n+1) &=& \mathcal{Y}(n) -\frac{1}{2}P_k(n)^{-1}\nabla F_k(n)
\end{eqnarray}
where $P_k(n) \coloneqq diag\{p_k(y_1(n)),\hdots,p_k(y_K(n))\}$. Therefore the iteration is of the form $\mathcal{Y}(n+1) = \mathcal{Y}(n) + \alpha_kd_k(n)$, where the descent direction $d_k(n) = -P_k(n)^{-1}\nabla F_k(n)$ satisfies $d_k(n)^T\nabla F_k(n)\leq 0$ with the equality being true only when $\nabla F_k(n) = 0$, which implies that the current iteration scheme is essentially a Descent Method. This observation allows us to analyze convergence rate of the proposed iteration scheme.

Let us now consider the scalability of the DA algorithm. If we scale the variables $\mathcal{Y}=\{y_j,j\in\mathbb{N}_{[1,K]}\}$ and $\mathcal{X}=\{x_i,i\in\mathbb{N}_{[1,N]}\}$ by a scaling factor $\sigma$, i.e., define $\mathcal{\tilde{Y}}=\{y_j/\sigma,j\in\mathbb{N}_{[1,K]}\}$ and $\mathcal{\tilde{X}}=\{x_i/\sigma,i\in\mathbb{N}_{[1,N]}\}$, the nature of optimization problem remains unchanged. In fact, minimizing the free-energy $F(\beta_k,\mathcal{X},\mathcal{Y})$ is equivalent to minimizing the modified free-energy $F(\beta_k\sigma^2,\mathcal{\tilde{X}},\mathcal{\tilde{Y}})$, and therefore the new iteration scheme is given by
\begin{equation}
	\mathcal{Y}(n+1) = \mathcal{Y}(n) - \frac{\sigma^2}{2}P_k(n)^{-1}\nabla F_k(n),
\end{equation}
where $\sigma$ can be appropriately chosen to obtain faster convergence rates. Another useful observation is that the scaling law relates the {\em spatial} scaling factor $\sigma$ to the annealing parameter $\beta$. This is intuitive, since to resolve smaller clusters one would require higher values of the annealing parameter $\tilde{\beta}_k = \beta_k\sigma^2$ to start with.

\section{CONCLUSIONS}\label{sec:Conclusion}

In this paper, we have discussed appropriate modifications of the Deterministic Annealing (DA) algorithm to address various heterogeneous capacity constraints. These modifications are achieved by appropriately modifying the free-energy function and the corresponding Gibbs distribution. The DA algorithm is found to achieve global minima (when verifiable) in several test simulations where other heuristics such as Llyod's algorithm fare poorly. We demonstrate the effectiveness of the proposed approach for some real-world scenarios and comparison against the previous heuristics reflect the strength of the proposed algorithm. We also establish equivalence to Descent Method and identify spatial scaling laws for faster convergence.

\addtolength{\textheight}{-12cm}   





\section*{ACKNOWLEDGMENT}
The authors would like to acknowledge NSF grants ECCS 15-09302, CMMI 14-63239 and CNS 15-44635 for supporting this work.


\bibliographystyle{IEEEtran} 
\bibliography{myRef}

\end{document}